\documentclass[8pt,a4paper]{article}
\usepackage{amsmath,amsfonts,amsthm,amsopn,color,amssymb}
\newcommand{\eps}{\varepsilon}

\newcommand{\R}{\mathbb{R}}
\newcommand{\Rn}{{\mathbb{R}^n}}

\newcommand{\Rnu}{{\mathbb{R}^{n-1}}}

\renewcommand{\a }{\alpha }
\renewcommand{\le }{\leqslant }
\renewcommand{\ge }{\geqslant }

\renewcommand{\O}{\Omega}
\newcommand{\G}{\Gamma}

\renewcommand{\i}{{H}^{1}(\Omega)}

\newcommand{\io}{\int _\O}
\newcommand{\iu}{\int _{\G_1}}
\newcommand{\id}{\int_{\G_2}}
\newcommand{\siD}{H^{\frac{1}{2}}_{00}(D)}
\newcommand{\siDd}{H^{\frac{1}{2}}_{00}(D)^*}
\newcommand{\sid}{H^{\frac{1}{2}}_{00}(\Gamma_2)}
\newcommand{\siu}{H^{\frac{1}{2}}_{00}(\Gamma_1)}
\newcommand{\sidd}{H^{\frac{1}{2}}_{00}(\Gamma_2)^*}
\newcommand{\sis}{H^{\frac{1}{2}}(\Sigma)}
\newcommand{\sisd}{H^{\frac{1}{2}}(\Sigma)^*}
\newcommand{\sig}{H^{\frac{1}{2}}(\Gamma)}
\newcommand{\siud}{H^{\frac{1}{2}}_{00}(\Gamma_1)^*}

\newcommand{\sidrd}{H^{\frac{1}{2}}_{00}(\G_{2,\rho})^*}
\newcommand{\si}{H^{\frac{1}{2}}(\partial \O)}
\newcommand{\n}{\nabla}

\newtheorem{theorem}{Theorem}[section]
\newtheorem{lemma}[theorem]{Lemma}

\newtheorem{corollary}[theorem]{Corollary}

\renewenvironment{proof}{\noindent{\textbf{Proof.\quad}}}{$\hfill\square$\vspace{0.2 cm}\\}

\begin{document}
{\title{ Solving elliptic Cauchy problems and the identification of nonlinear corrosion \thanks{Work supported in part by MIUR, grant n.
2002013279.} }}
\author{G. Alessandrini\thanks{Dipartimento di  Matematica e Informatica, Universit\`a degli Studi di Trieste, Italy,
 \tt{alessang@univ.trieste.it}} \ and \ E.
Sincich\thanks{S.I.S.S.A., Trieste, Italy, \tt { sincich@sissa.it}
}}

\date{}

\maketitle

\begin{abstract}
We deal with an inverse problem arising in corrosion detection. The presence of corrosion damage is modeled by a nonlinear boundary condition on the inaccessible portion of the metal specimen. We propose a method for the approximate reconstruction of such a nonlinearity. A crucial step of this procedure, which encapsulates the major cause of the ill-posedness of the problem, consists of the solution of a Cauchy problem for an elliptic equation. For this purpose we propose an SVD approach.

 \end{abstract}



\section{Introduction}
In this paper we deal with an inverse problem originating from corrosion detection. The corresponding direct problem, which models the electrochemical phenomenon of corrosion, is given as follows
\begin{equation}\label{P}
\left\{
\begin{array}
{lcl}
\Delta u=0\ ,& \mbox{in $\Omega$ ,}
\\
\dfrac{\partial u}{\partial \nu}=g\ ,   & \mbox{on $\G_2$ ,}
\\
\dfrac{\partial u}{\partial \nu}=f(u)\ , & \mbox{on $\G_1$ ,}
\\
u=0\ , & \mbox{on $\G_{D}$ ,}
\end{array}
\right.
\end{equation}
where $\G_1$ and $\G_2$
are two open, disjoint portions of  $\partial\O$ such that \mbox{${\G_{D}\!=\partial \!\O \setminus \!(\G_1\!\cup\!\G_2)}$}.
 We recall that, given $g\in L^2(\G_2)$, a weak solution of problem \eqref{P} is a function $u\in\i$, such that $u|_{\G_D}=0$ in the trace sense and which satisfies 
\begin{eqnarray}\label{soluzione}
\io \n u\cdot \n \eta  =\id g  \eta + \iu f(u) \eta \ ,
\end{eqnarray}
for every $\eta\in \i$ such that $\eta|_{\G_D}=0$.

 In the above boundary value problem $\O$ represents the electrostatic conductor, $\G_1$ represents the part of the boundary subject to corrosion, $\G_2$ represents the portion of the boundary accessible to direct inspection and $\G_D$ is a portion of the boundary where the electrostatic potential $u$ is grounded.
Such type of model with a specific choice of the nonlinear boundary term on $\G_1$ has been introduced and discussed by M.Vogelius and others \cite{Vog1},\cite{Vog2},\cite{Vog3}.

The inverse problem that we want to address here is the determination of the nonlinearity $f=f(u)$, when one non trivial pair of Cauchy data $u|_{\G_2},\frac{\partial u}{\partial \nu}|_{\G_2}$ is available. That is, we assume that we can measure, on the accessible boundary $\G_2$, the voltage and the current density. Note that also such an inverse problem must be formulated in a weak sense and the Cauchy data $\psi,\ g$ must be taken in the appropriate trace spaces, see Section 4 for details. In a previous paper \cite{as}, we have considered the issue of stability. The main results obtained there are summarized in the following Section 3. Let us just mention here, that we have considered classes of unknown nonlinearities for which also the direct problem \eqref{P} might not be well-posed, and also that for such an inverse problem also the domain within $\R$ where $f$ can be identified is part of the unknowns.

 In this paper we intend to initiate the study of a procedure for the approximate identification of the nonlinearity from approximate measurements of the data $u|_{\G_2},\frac{\partial u}{\partial \nu}|_{\G_2}$.  It seems necessary, as a first step of such a procedure, to solve the Cauchy problem for $u$ with Cauchy data on $\G_2$ and determine the corresponding Cauchy data for $u$ on the inaccessible part of the boundary $\G_1$, where the corrosion takes place. The approximate solution of a Cauchy problem for elliptic equations has already been studied by many authors, just to mention some of the most recent contributions, \cite{BE}, \cite{CHW}, \cite{EL}, \cite{HL}, \cite{KK}, \cite{KM}, \cite{KMF}, \cite{Le}, \cite{ML}, \cite{MEHILW}. 
In this paper we wish to propose an approach based on the reduction of the Cauchy problem to the regularized inversion of a suitable compact operator by the use of its singular value decomposition (SVD). In Section 4 we discuss this approach treating the problem in the wider generality of a Cauchy problems for variable coefficients elliptic equations. 
In Section 5 we specialize the same approach to the Laplace equation in a domain with a cylindrical geometry which might be well-suited to a reference conductor specimen, and to the model of electrochemical corrosion.
 We conclude in Section 6 by outlining the various steps of the approximate identification of the nonlinearity $f$.

\section{Main assumptions}

{\it A priori information on the domain}
\\
\\
We shall assume throughout that $\O$ is a bounded, connected domain in $\Rn$, $n\ge2$ such that $\mbox{diam} (\O)\le
D$ and with Lipschitz boundary $\partial \O$ with constants $r_0,
M$. 
Moreover, we assume that the portions of the boundary $\Gamma_i$ are contained respectively into surfaces $S_i$, $i=1,2$ which are $C^{1,\a}$ smooth with constants $r_0, M$.

We also suppose that the boundary of $\G_i$, within $S_i$, is of $C^{1,\a}$ class with constants $r_0, M$.

We introduce some notation that we shall use in the sequel, for every $\rho>0$ and $i=1,2$, we set
\begin{eqnarray}\label{ui}
U^i_{\rho}=\{x \in \bar{\O}:\mbox{dist}(x,\partial\O\setminus \G_i)>\rho \}\ \ ,
\end{eqnarray}
\begin{eqnarray}\label{gi}
\Gamma_{i,\rho}=U^1_{\rho}\cap \G_i\ .
\end{eqnarray} 
In some places, it will be necessary to isolate one privileged coordinate direction, to this purpose, we shall use the following notation for a point $x \in \Rn$,\\
 $x=(x',x_n)$, with $x'  \in \Rnu , \ x_n \in \mathbb{R}$.

 {\it A priori bound  on the energy}
\\
\\
We assume the following bound on the measured electrostatic potential 
\begin{eqnarray}\label{E}
\int_{\O}|\nabla u(x)|^2\le E^2\ .
\end{eqnarray}

 {\it A priori information on the boundary data}
\\
\\
The current flux $g$
is a prescribed function such that
\begin{eqnarray*}
 \|g\|_{C^{0,\a}(\G_2)}\le G\ , 
\end{eqnarray*}

 \begin{equation}\label{lower b}
\|g\|_{L^{\infty}(\G_{2,2r_0})}\ge m>0 \ .
\end{equation}

 {\it A priori information on the nonlinear term}
\\
\\
We assume that the function $f$ belongs to $C^{0,1}(\mathbb{R},\mathbb{R})$ and, in particular,
\begin{eqnarray}\label{f}
  f(0)=0\ \  \mbox{and}\ \  | f(u)-f(v)|\le L|u-v|\ \ \
 \mbox{for every } u, v\in \mathbb{R}\ .
\end{eqnarray}

{\it A priori information on the conductivity}
\\ 
\\ In Section 4 below, we shall consider Cauchy problems for solutions to variable coefficients elliptic equations of the form $\mbox{div}(\sigma\nabla u)=0\ \mbox{in}\ \O$. 
We shall assume that the conductivity tensor $\sigma(x)=(\sigma_{ij}(x))_{i,j=1}^{n}$ satisfies the ellipticity condition 
\begin{eqnarray}\label{lambda}
{\lambda}^{-1}|\xi|^2\le\sum_{i,j=1}^{n}\sigma_{ij}(x)\xi_i\xi_j\le \lambda|\xi|^2,\ \ \mbox{for all}\ \ x\in\O\ \ \mbox{and}\ \ \xi\in \Rn ,
\end{eqnarray}
 and the Lipschitz condition 
\begin{eqnarray}\label{Lipschitz}
|\sigma_{ij}(x)-\sigma_{ij}(y)|\le K|x-y|,\ \ \mbox{for all}\ i,j=1,\dots,n\  \mbox{and}\ x,y\in {\Omega},
\end{eqnarray}
 where $K>0, \lambda\ge1$ are prescribed constants.

From now on we shall refer to the \emph{a priori data} as to the set of quantities
$r_0, M, \a, L, G, E, D, m, \lambda, K$.\\
In the sequel we shall denote by $\eta(t)$ and $\omega(t)$, two positive increasing functions defined on $(0, +\infty)$, that satisfy
\begin{eqnarray}\label{eta}
&&\eta(t)\ge\exp{\bigg[-\bigg(\frac{t}{c}\bigg)^{-\gamma}\bigg]}, \ \ \ \mbox{for every}\ \ 0<t\le G\ \ ,
\end{eqnarray}
\begin{eqnarray}\label{omega}
&&\omega(t)\le C \left|\log (t) \right|^{-\theta},\ \ \ \mbox{for every}\ \ 0<t<1\ \ ,
\end{eqnarray}
where $c>0$, $C>0$, $\gamma>1$, $0<\theta<1$ are constants depending on the \emph{a priori data} only.
\section{Stability result}
We review here the stability result obtained in \cite{as}.
As preliminary step, we evaluate the amplitude of the range of $u$ on $\G_1$ and we prove a lower bound on the oscillation of $u$ on $\G_1$.
\begin{theorem}[Lower bound for the oscillation]\label{teo1}
Let $\O, g $ satisfying the a priori assumptions. Let $u$ be a
weak solution of \eqref{P} satisfying the a priori bound \eqref{E}
then $$\mathop{\rm osc}\limits_{\G_1} u\ge
\eta(\|g\|_{L^{\infty}(\G_{2,2r_0})})$$ where $\eta$ satisfies
\eqref{eta}.
\end{theorem}
\begin{proof}
[Sketch]. First we note that the (unknown) Cauchy data $u|_{\G_1}, {\frac{\partial u}{\partial \nu}}|_{\G_1}$ can be dominated by the oscillation of $u$ on $\G_1$. Next, we use a stability result for the Cauchy problem with data on $\G_1$ and obtain a bound on ${\frac{\partial u}{\partial \nu}}|_{\G_2}$. We refer to \cite[Section 3]{as}.
\end{proof}
The main result in \cite{as} is the following.
\begin{theorem}[Stability for the nonlinear term $f$]\label{teo2}
Let \mbox{$u_i \in \i, \ i=1,2$}\ be weak solutions of the problem $\eqref{P}$, with $f=f_i$ and $g=g_i$
respectively and such that \eqref{E} holds for each $u_i$.
Let us also assume that, for some positive number $m$, the following holds
\begin{eqnarray}\label{m}
\|g_1\|_{L^{\infty}(\G_{2,2r_0})}\ge m>0\ .
\end{eqnarray}
Moreover, let $\psi_i=u_i\big|_{\G_2}, \ i=1,2$.
There exists $\eps_{0}>0$ only depending on the \emph{a priori data} and on $m$ such that, if, for some $\eps,\ 0<\eps<\eps_0$, we have
\begin{eqnarray*}
&& \|\psi_1 -\psi_2  \|_{L^{2}(\G_2)}\le \eps\ , \\
&& \|g_1 -g_2  \|_{L^{2}(\G_2)}\le \eps\ ,
\end{eqnarray*}
 then
$$\|f_1-f_2 \|_{L^{\infty}(V)}\le \omega(\eps)\ , $$
where
 $$V=(\a,\beta)\subseteq [-CE,CE]\ ,$$
is  such that
 $$\beta -\a >\frac{\eta(m)}{2}$$
and $\eta, \omega$ satisfy \eqref{eta}, \eqref{omega} respectively.
\end{theorem}
\begin{proof}
[Sketch]. In this case we use first a stability estimate for a Cauchy problem with data on $\G_2$. Next we show that on a suitable curve on $\G_1$, the restriction of $u_1, u_2$ are strictly monotone and their ranges agree on a sufficiently large interval $V$. On such an interval we are able to bound $f_1-f_2$ in terms of $(u_1-u_2)|_{\G_1}$ and $({\frac{\partial u_1}{\partial \nu}}-{\frac{\partial u_2}{\partial \nu}})|_{\G_1}$. See \cite[Section 4]{as} for details.
\end{proof}

\section{Solving the Cauchy problem}

We consider here a Cauchy problem for an elliptic equation with variable coefficients
\begin{equation}\label{cpg}
\left\{
\begin{array}
{lcl}
\mbox{div}(\sigma \nabla u)=0& \mbox{in $\Omega$},
\\
u=\psi  & \mbox{on $\Sigma$},
\\
\sigma\nabla u\cdot\nu=g & \mbox{on $\Sigma$ .}
\end{array}
\right.
\end{equation}
Here the conductivity tensor $\sigma(x)=(\sigma_{ij}(x))_{i,j=1}^{n}$ satisfies the ellipticity condition \eqref{lambda} and the Lipschitz condition \eqref{Lipschitz}. The domain $\O$ satisfies the same assumptions stated previously in Section 2 and $\Sigma$ is an open connected portion of $\partial \O$ which is $C^{1,\a}$ smooth with constant $r_0, M$ as it was previously stated for the portions $\G_1, \G_2$. We also suppose that the boundary of $\Sigma$ within $\partial \O$ is of Lipschitz class with constants $r_0, M$. Moreover, we denote with $U^{\Sigma}_{\rho}$ and $\Sigma_{\rho}$ the analogous of the sets defined in \eqref{ui} and \eqref{gi} respectively, with $\G_i$ replaced by $\Sigma$.
We introduce the trace spaces $\sis, \ H^{\frac{1}{2}}_{00}(\Sigma)$ as the interpolation spaces $[H^1(\Sigma),L^2(\Sigma)]_{\frac{1}{2}}$, $[H_{0}^1(\Sigma),L^2(\Sigma)]_{\frac{1}{2}}$ respectively, see \cite[Chap. 1]{LiMa} for details.
We shall denote the corresponding dual spaces by $\sisd, \ H_{00}^{\frac{1}{2}}(\Sigma)^*$, respectively.

 We recall that there exists a linear extension operator
\begin{eqnarray}\label{esteE}
&&E:\sis \rightarrow \si\ ,\ \ \mbox{such that} \ \ E(\psi)=\psi\ \ \ \mbox{on}\ \Sigma\ \ \mbox{and}\nonumber\\
&&\|E(\psi)\|_{\si}\le C\|\psi\|_{\sis}\, \mbox{for every}\ \psi\in\sis,
\end{eqnarray}
where $C>0$ is a constant depending on the \emph{a priori data} only, see for instance \cite[Lemma 7.45]{Adams}. Also we recall that the operator $E_0$ of continuation to zero outside $\Sigma$,
\begin{equation}\label{estensioneazero}
E_0(\psi)=\left\{
\begin{array}{rl}
\psi,\ & \mbox{in}\ \Sigma, \\
0,\  & \mbox{in}\ \partial\O\setminus\Sigma,
\end{array}
\right.
\end{equation}
is bounded from $H^{\frac{1}{2}}_{00}(\Sigma)$ into $\si$. Note that, by such an extension, $H_{00}^{\frac{1}{2}}(\Sigma)$ can be identified with the closed subspace of $\si$ of functions supported in $\overline{\Sigma}\subset\partial\O$. More precisely, denoting by $\Gamma=\partial\O\setminus\overline{\Sigma}$ and 
\begin{eqnarray}\label{omegagamma}
H^{1}_0(\O,\G)=\{u\in H^1(\O):\ u|_{\G}=0\ \mbox{in the trace sense} \}
\end{eqnarray}
we can identify $H^{\frac{1}{2}}_{00}(\Sigma)$ with the trace space of $H^{1}_0(\O,\G)$ on $\partial\O$. See \cite[Chap. 1]{LiMa} and also, for more details, \cite{Tar}.

Given $\psi\in\sis$ and $g\in H^{\frac{1}{2}}_{00}(\Sigma)^*$ we shall say that $u\in H^1(\O)$ is a weak solution to \eqref{cpg} if $u|_{\Sigma}=\psi$ in the trace sense and also 
\begin{eqnarray}
\int_{\O}\sigma\nabla u\cdot\nabla \eta=<g,\eta|_{\Sigma}>
\end{eqnarray}
for every $\eta\in H^1_0(\O,\Gamma)$. Here $<\cdot,\cdot>$ denotes the pairing between $H^{\frac{1}{2}}_{00}(\Sigma)^*$ and $H^{\frac{1}{2}}_{00}(\Sigma)$ based on the $L^2(\Sigma)$ scalar product.
Our first step in the solution of the Cauchy problem \eqref{cpg} is the reduction to the case when $\psi=0$.
To this purpose we consider the weak solution $W\in H^1(\O)$ to the well-posed Dirichlet problem 
\begin{equation}\label{Dirichlet}
\left\{
\begin{array}
{lcl}
\mbox{div}(\sigma \nabla W)=0& \mbox{in $\Omega$},
\\
W=E{\psi}  & \mbox{on $\partial \O$}.
\end{array}
\right.
\end{equation}
Setting $U=u-W$ and $G=g-\sigma\nabla W\cdot\nu|_{H^{\frac{1}{2}}_{00}(\Sigma)}\in H^{\frac{1}{2}}_{00}(\Sigma)^*$, we have that $U$ is a weak solution to the Cauchy problem 
\begin{equation}\label{omogeneog}
\left\{
\begin{array}
{lcl}
\mbox{div}(\sigma \nabla U)=0& \mbox{in $\Omega$},
\\
U=0  & \mbox{on $\Sigma$},
\\
\sigma\nabla U\cdot\nu=G & \mbox{on $\Sigma$} .
\end{array}
\right.
\end{equation}
For every $h\in H^{\frac{1}{2}}_{00}(\Gamma)^*$ let us consider the mixed boundary value problem 
\begin{equation}\label{compattog}
\left\{
\begin{array}
{lcl}
\mbox{div}(\sigma \nabla v)=0& \mbox{in $\Omega$},
\\
v=0  & \mbox{on $\Sigma$},
\\
\sigma\nabla v\cdot\nu=h& \mbox{on $\Gamma$} \ .
\end{array}
\right.
\end{equation}
Here, denoting by $H^1_0(\Omega,\Sigma)$ the space introduced in \eqref{omegagamma} when $\G$ replaced with $\Sigma$, a function $v\in H^1_0(\O,\Sigma)$ is said to be a weak solution to \eqref{compattog} if
\begin{eqnarray}\label{formulazione1}
\int_{\O}\sigma\nabla v\cdot \nabla \eta=<h,\eta|_{\G}>\ \ \mbox{for every}\ \eta\in H^1_0(\O,\Sigma) .
\end{eqnarray}
It is readily seen, by the Lax-Milgram Theorem, that such mixed boundary value problem \eqref{compattog} is well-posed. It is also evident that, finding the appropriate $h \in H^{\frac{1}{2}}_{00}(\Gamma)^*$ such that $\sigma\nabla v\cdot\nu|_{H^{\frac{1}{2}}_{00}(\Gamma)}=G$, would imply that $v=U$ and provide us with the solution to \eqref{omogeneog}. We note however, that given $\rho_0>0$ such that $\Sigma_{\rho_0}$ has nonempty interior, it would suffice to check that for some $\rho, 0<\rho<\rho_0$, $\sigma\nabla v\cdot\nu=G$ when both functionals are restricted to $H^{\frac{1}{2}}_{00}(\Sigma_{\rho}).$ In fact, this is a consequence of the uniqueness of the solution of the Cauchy problem when the Cauchy data are prescribed on $\Sigma_{\rho}$ (instead than on all of $\Sigma$). Thus, having fixed $\rho,\ 0<\rho<\rho_0$, the solution of the Cauchy problem \eqref{omogeneog} amounts to find $h\in H^{\frac{1}{2}}_{00}(\Gamma)^*$ such that $\sigma\nabla v\cdot\nu=G$ on $H^{\frac{1}{2}}_{00}(\Sigma_{\rho})$.
   
We prove the following.
\begin{theorem}\label{operatorecompatto}
For any $\rho, 0<\rho<\rho_0$, let $T_{\rho}$ be the operator 
\begin{eqnarray}\label{opc}
T_{\rho}:H^{\frac{1}{2}}_{00}(\Gamma)^* &\rightarrow & H^{\frac{1}{2}}_{00}(\Sigma_{\rho})^* \\
 h &\mapsto  & \left.\sigma\nabla v\cdot\nu\displaystyle \right|_{\Sigma_{\rho}}\nonumber
\end{eqnarray}
where $v\in H^1_0(\O,\Sigma)$ solves the mixed problem \eqref{compattog}. The operator $T_{\rho}$ is compact.
\end{theorem} 
\begin{proof}
By the well posedness of the mixed boundary value problem \eqref{compattog}, 
 the linear operator 
$$
\begin{array}{ccc}
S: H^{\frac{1}{2}}_{00}(\Gamma)^* &\rightarrow & H^{1}_0(\O,\Sigma) \\
 h &\mapsto  & v
\end{array}
$$
is bounded.

Moreover, by a standard result of regularity at the boundary, it follows that for every $\rho>0,\ v\in C^{1,\a}(U_{\rho}^{\Sigma})$ and there exists a constant $C_{\rho}>0$ depending on the \emph{a priori data} and on $\rho$ only, such that 
$$\|v\|_{{C^{1,\a}}(\Sigma_{\rho})}\le C_{\rho}\|v\|_{H^1_0(\O)}.$$  
Thus the operator
$$
\begin{array}{ccc}
D_{\rho}:\i\ &\rightarrow & C^{0,\a}(\Sigma_{\rho})  \\
 v &\mapsto  & \left.\sigma\nabla v\cdot\nu\displaystyle \right|_{\Sigma_{\rho}}
\end{array}
$$
is bounded.
Finally, since the inclusion 
$$
\begin{array}{ccc}
i_{\rho}:\ C^{0,\a}(\Sigma_{\rho})\hookrightarrow & H^{\frac{1}{2}}_{00}(\Sigma_{\rho})^* \\
 \end{array}
$$
is compact and $T_{\rho}$ can be factored as 
$T_{\rho}=i_{\rho}\circ D_{\rho}\circ S $, the thesis follows.
\end{proof}
Being $T_{\rho}$ a compact operator between Hilbert spaces, we have that there exists a triple $\{{\sigma}_j^{\rho}, h_j, {g}_j^{\rho} \}_{j=1}^{\infty}$ called \emph{singular value decomposition}, such that ${\{\sigma_j^{\rho}\}}_{j=1}^{\infty}$ is a non increasing infinitesimal sequence of nonnegative numbers,  ${\{h_j\}}_{j=1}^{\infty}$, ${\{{g}_j^{\rho}\}}_{j=1}^{\infty}$ are orthonormal bases for $H^{\frac{1}{2}}_{00}(\Gamma)^*$ and for $ H^{\frac{1}{2}}_{00}(\Sigma_{\rho})^*$ respectively, and moreover it holds
\begin{eqnarray}\label{scomposizione}
T_{\rho}h_j={\sigma}_j^{\rho}{g}_j^{\rho},\ \ \mbox{for every}\ \ j=1,2,\dots \\
T^*_{\rho}{g}_j^{\rho}={\sigma}_j^{\rho}h_j,\ \ \mbox{for every}\ \ j=1,2,\dots 
\end{eqnarray}
where $T^*_{\rho}$ denotes the adjoint operator to $T_{\rho}$.
 By the regularization theory for the inversion of compact operators, we have that, denoting with $(\cdot,\cdot)_{H^{\frac{1}{2}}_{00}(\Sigma_{\rho})^*}$ the scalar product for the Hilbert space $H^{\frac{1}{2}}_{00}(\Sigma_{\rho})^*$, the family of operators $R_{\a},\ \a>0$
\begin{eqnarray}\label{ra}
\begin{array}{ccc}
R_{\a}: H^{\frac{1}{2}}_{00}(\Sigma_{\rho})^*&\rightarrow &H^{\frac{1}{2}}_{00}(\Gamma)^*  \\
 g &\mapsto  & \sum_{{\sigma}_k^{\rho}\ge \a}\frac{1}{{\sigma}_k^{\rho}}{(g,{g}_k^{\rho})}_{H^{\frac{1}{2}}_{00}(\Sigma_{\rho})^*}h_k
\end{array}
\end{eqnarray}
is a regularization strategy for $T_{\rho}$, namely
\begin{eqnarray}\label{lim}
\lim_{\a\rightarrow 0} R_{\a}T_{\rho}h=h\ ,\ \mbox{for every}\ h\in H^{\frac{1}{2}}_{00}(\Gamma)^* ,
\end{eqnarray}
(see for instance \cite [Chap. 2]{Kir}).
Moreover, the choice 
\begin{eqnarray}\label{par}
\a(\eps)=\eps^{2(1-\gamma)}
\end{eqnarray}
 where $\gamma$ is a fixed number, $0<\gamma<1$, is an admissible one, this means that if given, for every $\eps>0$, $g, g_{\eps}\in H^{\frac{1}{2}}_{00}(\Sigma_{\rho})^*$ and $h \in H^{\frac{1}{2}}_{00}(\Gamma)^*$ such that 
\begin{eqnarray}
g=T_{\rho}h\ \ \ \mbox{and}\ \ \ \|g-g_{\eps}\|_{H^{\frac{1}{2}}_{00}(\Sigma_{\rho})^*}\le \eps\ ,
\end{eqnarray}
then it follows that 
\begin{eqnarray}\label{reg}
\lim_{\eps\rightarrow 0}\|R_{\alpha(\eps)}g_{\eps}-h\|_{H^{\frac{1}{2}}_{00}(\Gamma)^*}=0\ .
\end{eqnarray}
We can return now to the Cauchy problem \eqref{cpg}, when $\psi$ is arbitrary in $H^{\frac{1}{2}}(\Sigma)$.
Let us suppose that, for every $\eps>0$, $\psi_{\eps}\in \sis$, $g_{\eps}\in H^{\frac{1}{2}}_{00}(\Sigma_{\rho})^*$, and let $W_{\eps}\in \i$ be the weak solution of \eqref{Dirichlet}, with ${\psi}={\psi_{\eps}}$.
Let us denote by $R^{\eps}=R_{\a(\eps)}(g_{\eps}-\left.\sigma\nabla W_{\eps}\cdot\nu\right|_{\Sigma_{\rho}})+ \left.\sigma\nabla W_{\eps}\cdot\nu\right|_{\Gamma}\in H^{\frac{1}{2}}_{00}(\Gamma)^*$, where $R_{\a}$ and $\a(\eps)$ are the regularization strategy and the regularization parameter introduced in \eqref{ra} and \eqref{par}, respectively. 
We propose as approximate regularized solution to the problem \eqref{cpg} the function $u_{\eps}\in \i$ which is a weak solution of the mixed boundary value problem
\begin{equation}\label{mbvp}
\left\{
\begin{array}
{lll}
\mbox{div}(\sigma\nabla u_{\eps})=0& \mbox{in $\Omega$},
\\
u_{\eps}=\psi_{\eps}  & \mbox{on $\Sigma$},
\\
\sigma\nabla u_{\eps}\cdot\nu=R^{\eps}& \mbox{on $\Gamma$.}
\end{array}
\right.
\end{equation}
In analogy to \eqref{compattog} and \eqref{formulazione1}, we shall call weak solution of the problem \eqref{mbvp}, a function $u_{\eps}\in H^1(\O)$ such that $u_{\eps}|_{\Sigma}=\psi_{\eps}$ in the trace sense and such that 
\begin{eqnarray}\label{formulazione2}
\int_{\O}\sigma\nabla u_{\eps}\cdot\nabla \eta=<R^{\eps}, \eta|_{\Gamma}>\ \ \mbox{for every}\ \  \eta \in H^1_0(\O,\Sigma).
\end{eqnarray}
 The well-posedeness of problem \eqref{mbvp} is again a consequence of the Lax-Milgram Theorem. The following Theorem provides a convergence results for the procedure of regularized inversion of the Cauchy problem \eqref{cpg} that we have just outlined.
\begin{theorem}\label{approssimazione}
Let $\psi\in\sis$ and $g\in\ H^{\frac{1}{2}}_{00}(\Sigma)^*$ and suppose that there exists $u \in H^1(\O)$, which is a weak solution to the Cauchy problem \eqref{cpg}. If, given $\eps>0$, we have that $\psi_{\eps}\in\sis$ and $g_{\eps}\in H^{\frac{1}{2}}_{00}(\Sigma_{\rho})^*$
\begin{eqnarray}
&&\|\psi-\psi_{\eps}\|_{\sis}\le \eps\ , \label{dpb}\\
&&\|g-g_{\eps}\|_{H^{\frac{1}{2}}_{00}(\Sigma_{\rho})^*}\le \eps\ , \label{npb}
\end{eqnarray}
then
\begin{eqnarray}
&&\lim_{\eps\rightarrow 0}u_{\eps}|_{\Gamma}=u|_{\Gamma}\ \ \ in\ \ \sig\ \ \label{dpbt}\ ,\\
&&\lim_{\eps\rightarrow 0}\left.\sigma\nabla u_{\eps}\cdot\nu\right|_{\Gamma}=\left.\sigma\nabla u\cdot\nu\right|_{\Gamma}\ \ \ in\ \ {{H^{\frac{1}{2}}_{00}(\G)}^*} \label{npbt}\ .
\end{eqnarray}
\end{theorem} 
\begin{proof}
Let us observe that given $S$ any open and connected portion of $\partial \O$, the following holds
\begin{eqnarray*}
\left\|\left.\sigma\nabla W_{\eps}\cdot\nu\right|_{S}-\left.\sigma\nabla W\cdot\nu\right|_{S}\right\|_{H_{00}^{\frac{1}{2}}(S)^*}
\le c_1 \|W-W_{\eps}\|_{\i}\le c_2\|E{\psi_{\eps}}-E{\psi}\|_{\si}
\end{eqnarray*}
then replacing in \eqref{esteE} ${\psi}$ with ${\psi_{\eps}}-{\psi}$, we have by \eqref{dpb} that 
\begin{eqnarray}\label{bpW}
\left\|\left.\sigma\nabla W_{\eps}\cdot\nu\right|_{S}-\left.\sigma\nabla W\cdot\nu\right|_{S}\right\|_{H_{00}^{\frac{1}{2}}(S)^*}
\le c_3 \eps \ , 
\end{eqnarray}
where $c_1,c_2,c_3>0$ are constants depending on the \emph{a priori data} and on $S$ only.
Thus by \eqref{bpW}, with $S=\Sigma_{\rho}$, and by \eqref{npb}, we have that 
\begin{eqnarray}\label{bpWg}
\lim_{\eps\rightarrow 0}\left\|g-g_{\eps}+\left.\sigma\nabla W_{\eps}\cdot\nu\right|_{\Sigma_{\rho}}-\left.\sigma\nabla W\cdot\nu\right|_{\Sigma_{\rho}}\right\|_{H_{00}^{\frac{1}{2}}(\Sigma_{\rho})^*}=0
\end{eqnarray}
Moreover, we have that \eqref{npbt} follows by applying \eqref{bpW} with $S=\Gamma$, \eqref{reg} with $g_{\eps}$ replaced with $g_{\eps}-\sigma\nabla W_{\eps}\cdot \nu|_{\Sigma_{\rho}}$ and \eqref{bpWg}. Indeed, we have
\begin{eqnarray*}
&&
\left\|\left.\sigma\nabla u\cdot\nu\right|_{\Gamma}-\left.\sigma\nabla u_{\eps}\cdot\nu\right|_{\Gamma}\right\|_{H^{\frac{1}{2}}_{00}(\G)^*}
\le\\
&&\le\left\|R_{\a(\eps)}\left(g_{\eps}-\left.\sigma\nabla W_{\eps}\cdot\nu\right|_{\Sigma_{\rho}}\right)+\left.\sigma\nabla W\cdot\nu\right|_{\G}-\left.\sigma\nabla u\cdot \nu\right|_{\Gamma}\right\|_{H^{\frac{1}{2}}_{00}(\G)^*}+\\
&&\ \ \ +\left\|\left.\sigma\nabla W\cdot\nu\right|_{\Gamma}-\left.\sigma\nabla W_{\eps}\cdot\nu\right|_{\Gamma}\right\|_{H^{\frac{1}{2}}_{00}(\Gamma)^*}.
\end{eqnarray*} 

Finally, by a standard trace inequality 
\begin{eqnarray}
&&\|u|_{\Gamma}-u_{\eps}|_{\Gamma}\|_{\sig}\le c_4\|u-u_{\eps}\|_{\i}\le\nonumber\\
&&\le c_5\left(\left\|\left.\sigma\nabla u\cdot\nu\right|_{\Gamma}-\left.\sigma\nabla u_{\eps}\cdot\nu\right|_{\Gamma}\right\|_{H^{\frac{1}{2}}_{00}(\G)^*} +\|\psi-\psi_{\eps}\|_{\sis}\right)
\end{eqnarray}
where $c_4, c_5>0$ are constants depending on the \emph{a priori data} only,
then \eqref{dpbt} follows by recalling \eqref{npbt} and from \eqref{dpb}.
\end{proof}

\section {A special case}

Let $D$ be a bounded domain in $\mathbb{R}^{n-1}$, with Lipschitz boundary $\partial D$ with constants $r_0,M$. 
From now on we shall consider this special choice of $\O$
\begin{eqnarray*}
\O=D\times(0,1)\ ,\ \G_2=D\times\{0\}\ ,\ \G_1=D\times\{1\}\ ,\ \G_D=\partial D\times (0,1).
\end{eqnarray*}
In the following we will denote by $\lambda_k, \varphi_k,\ k=1,2,\dots$, the Dirichlet eigenvalues and eigenfunctions of $-\Delta$ on $D$, namely
\begin{equation}
\left\{
\begin{array}
{lll}
-\Delta \varphi_k=\lambda_k\varphi_k& \mbox{in $D$},
\\
\ \ \ \varphi_k \in H_{0}^{1}(D)\ .
\end{array}
\right.
\end{equation}

We recall that the family $\{\varphi_k\}_{k=1}^{\infty}$ is an orthogonal basis in $L^2(D)$ and also in $H^1_0(D)$. In the following we shall refer to the $\{\varphi_k\}_{k=1}^{\infty}$ as the basis normalized in the $L^2(D)$ norm.
 We have that $\psi\in H^{\frac{1}{2}}_{00}(D)$ if and only if its Fourier coefficients
\begin{eqnarray}
\psi_k=\int_{D}\psi\varphi_k
\end{eqnarray}
 satisfy 
\begin{eqnarray}
\sum_{k=1}^{\infty}{\lambda_k}^{\frac{1}{2}}{\psi_k}^2<\infty 
\end{eqnarray}
and that, as a norm on $H^{\frac{1}{2}}_{00}(D)$ we can choose 
\begin{eqnarray}
\|\psi \|_{H^{\frac{1}{2}}_{00}(D)}=\left(\sum_{k=1}^{\infty}{\lambda_k}^{\frac{1}{2}}{\psi_k}^2\right)^{\frac{1}{2}}\ .
\end{eqnarray}
Moreover, $h\in H^{\frac{1}{2}}_{00}(D)^*$ if and only if, its  Fourier coefficients
\begin{eqnarray}\label{pairing}
h_k=<h,\varphi_k>\ ,
\end{eqnarray}
 satisfy 
\begin{eqnarray}
\sum_{k=1}^{\infty}{\lambda_k}^{-\frac{1}{2}}{h_k}^2<\infty 
\end{eqnarray}
and the norm on $H^{\frac{1}{2}}_{00}(D)^*$ turns out to be  
\begin{eqnarray}
\|h\|_{H^{\frac{1}{2}}_{00}(D)^*}=\left(\sum_{k=1}^{\infty}{\lambda_k}^{-\frac{1}{2}}{h_k}^2\right)^{\frac{1}{2}}.
\end{eqnarray}
Here $<\cdot,\cdot>$ denotes the pairing between $H^{\frac{1}{2}}_{00}(D)^*$ and $H^{\frac{1}{2}}_{00}(D)$ based on the $L^2(D)$ scalar product.
Note also that $\{{\lambda_k}^{-\frac{1}{4}}\varphi_k\}$ and $\{{\lambda_k}^{\frac{1}{4}}\varphi_k\}$ constitute orthonormal bases for $H^{\frac{1}{2}}_{00}(D)$ and $H^{\frac{1}{2}}_{00}(D)^*$ respectively.

Due to the cylindrical geometry of $\O$, we remark that we can identify the spaces $H^{\frac{1}{2}}_{00}(\G_i),\ H^{\frac{1}{2}}_{00}(\G_i)^*,\ i=1,2$, with $H^{\frac{1}{2}}_{00}(D),\ H^{\frac{1}{2}}_{00}(D)^*$ respectively. Furthermore, as noted in Section 4, we can identify $H^{\frac{1}{2}}_{00}(\G_1)$ with the trace space on $\partial \O$ of $H^1_0(\O,\G)$ when $\G=\stackrel{\circ}{(\overline{\G_2\cup\G_D})}$, and the same holds when the roles of $\G_1$ and $\G_2$ are exchanged.

Let  $\psi\ \in \sid$, $g \in \sidd$ and let us consider the following Cauchy problem with auxiliary homogeneous condition on $\G_D$
\begin{equation}\label{cp}
\left\{
\begin{array}
{lcl}
\Delta u=0& \mbox{in $\Omega$},
\\
u=\psi  & \mbox{on $\G_2$},
\\

\dfrac{\partial u}{\partial \nu}=g & \mbox{on $\G_2$},
\\
u=0 & \mbox{on $\G_{D}$}.
\end{array}
\right.
\end{equation}
We shall say that $u$ is a weak solution to the problem \eqref{cp} if $u|_{\stackrel{\circ}{(\overline{\G_2\cup\G_D})}}=E_0(\psi)$ in the trace sense and if 
$$\int_{\O}\nabla u\cdot\nabla \eta=<g, \eta|_{\G_2}>\ \ \mbox{for every}\ \eta \in \ H^1_0(\O,\stackrel{\circ}{(\overline{\G_1\cup\G_D})}).$$
Here $E_0(\psi)$ denotes the extension of $\psi$ by zero outside $\G_2$ and $<\cdot,\cdot>$ denotes the pairing between $H^{\frac{1}{2}}_{00}(\G_2)^*$ and $H^{\frac{1}{2}}_{00}(\G_2)$ based on the $L^2(\G_2)$ scalar product. 
We shall use a strategy similar to the one discussed in Section 4, but with some slight variations, suggested by the presence of the portion $\G_D$ of the boundary where $u=0$. As before, we reduce the problem \eqref{cp} to the special case when $\psi=0$ and introduce the well-posed Dirichlet problem
\begin{equation}\label{1}
\left\{
\begin{array}{ll}
\Delta v=0& \mbox{in $\Omega$},
\\
v=\xi  & \mbox{on $\G_1$},
\\
v=0   & \mbox{on $\stackrel{\circ}{(\overline{\G_2\cup\G_D})},$}
\end{array}
\right.
\end{equation}
where $\xi$ is a prescribed function in  $H^{\frac{1}{2}}_{00}(\G_1)$. 
To this purpose, in analogy with \eqref{Dirichlet}, we consider $W\in\i$ as the weak solution to the Dirichlet problem
\begin{equation}\label{2}
\left\{
\begin{array}{ll}
\Delta W=0& \mbox{in $\Omega$},
\\
W=\psi  & \mbox{on $\G_2$},
\\
W=0   & \mbox{on $\stackrel{\circ}{(\overline{\G_1\cup\G_D})}.$}
\end{array}
\right.
\end{equation}
The difference $U=u-W$ shall satisfy \eqref{cp} with $\psi=0$ and $g$ replaced with $G=g-\frac{\partial W}{\partial \nu}|_{H^{\frac{1}{2}}_{00}(\G_2)}$.

Note that the well posed boundary value problem \eqref{1}, will take the place of \eqref{compattog}.
We intend to invert the map 
\begin{eqnarray}
T:\xi\rightarrow \left.\frac{\partial v}{\partial \nu}\right|_{\G_2}
\end{eqnarray}
in order to solve the Cauchy problem.
It is convenient at this stage to recall the identification of the trace spaces on $\G_i,\ i=1,2$ with the corresponding ones on $D$.

\begin{lemma}\label{compatto}
Let $T$ be the operator 

\begin{eqnarray}
T: \siD  &\rightarrow & \siDd \\
 \xi &\mapsto  & \left.\frac{\partial v}{\partial \nu}\displaystyle \right|_{\G_2}
\end{eqnarray}

where $v$ is the weak solution of the problem \eqref{1}.
Then $T$ extends to a compact and self-adjoint operator on $L^2(D)$, such that $\left\{-{\lambda_k}^{\frac{1}{2}}({\sinh ({\lambda_k}^{\frac{1}{2}})})^{-1},\varphi_k\right\}_{k=1}^{\infty}$ are its eigenvalues and eigenfunctions respectively. The singular value decomposition of $T: \siD  \rightarrow  \siDd$ is given by 
\begin{eqnarray}\label{svd}
 \{-({\sinh({\lambda_k}^{\frac{1}{2}})})^{-1},{\lambda_k}^{-\frac{1}{4}}\varphi_k,{\lambda_k}^{\frac{1}{4}}\varphi_k \}_{k=1}^{\infty}\ .
\end{eqnarray}
\end{lemma}
\begin{proof}
Let us first observe that the operator $T$ is well defined since the problem \eqref{1} is well-posed.
In this special setting we can represent the solution $v$ of \eqref{1} by separation of variables, namely
\begin{eqnarray}
v(x',x_n)=\sum_{k=1}^{\infty}\frac{\xi_k}{\sinh({\lambda_k}^{\frac{1}{2}})}\sinh({\lambda_k}^{\frac{1}{2}}x_n)\varphi_k(x')
\end{eqnarray}
where $\{\xi_k\}_{k=1}^{\infty}$ are the Fourier coefficients of $\xi$ with respect to the $L^2(D)$ basis $\{\varphi_k\}_{k=1}^{\infty}$.
After straightforward calculations we have that 
\begin{eqnarray}\label{auto}
T\left(\sum_{k=1}^{\infty}\xi_k \varphi_k\right)=\sum_{k=1}^{\infty}\left(-\frac{\xi_k {\lambda_k}^{\frac{1}{2}}}{\sinh({\lambda_k}^{\frac{1}{2}})}\right)\varphi_k 
\end{eqnarray}
thus the operator extends to a self-adjoint operator on $L^2(D)$ and since the eigenvalues are infinitesimal we conclude that $T$ is compact as an operator from $L^2(D)$ into $L^2(D)$. Moreover, since $\siD$ is continuously embedded in $L^2(D)$ and $L^2(D)$ is continuously embedded in $\siDd$, also $T: \siD  \rightarrow  \siDd$ is compact and its SVD turns out to be \eqref{svd}. 
\end{proof}

As a consequence of the above Lemma \ref{compatto}, we obtain that the family of operators
\begin{eqnarray*}
R_{\alpha}:\siDd\  \longrightarrow \siD,\ \ \mbox{such that}\nonumber 
\end{eqnarray*}
\begin{eqnarray}\label{Rd}
R_{\alpha}(G)=\sum_{\mu_k\ge \alpha}(-{\sinh ({\lambda_k}^{\frac{1}{2}})})(G,\varphi_k)_{\siDd}\varphi_k
\end{eqnarray}
where $\mu_k=({\sinh({\lambda_k}^{\frac{1}{2}})})^{-1}$, is a regularization strategy for $T$ and the choice \eqref{par} for the parameter $\a$ is still admissible.
We are in the position now to present the regularized approximate solution for the following special case of the problem \eqref{cp}. That is, given $G\in H^{\frac{1}{2}}_{00}(\G_2)$, 
\begin{equation}\label{3}
\left\{
\begin{array}{ll}
\Delta U=0& \mbox{in $\Omega$},
\\
U=0  & \mbox{on $\G_2$},
\\
\dfrac{\partial U}{\partial \nu}=G & \mbox{on $\G_2$},
\\
U=0   & \mbox{on $\G_D$.}
\end{array}
\right.
\end{equation}

In this section we shall denote by $[r]$ the integral part of the real number $r$.
\begin{theorem}\label{omogeneo}
For every $\eps>0$, let $G_{\eps}\in \sidd$ and suppose that there exists $U \in \i$, which is a weak solution of the problem \eqref{3}. If we have 
\begin{eqnarray*}
\|G_{\eps}-G\|_{\sidd}\le \eps
\end{eqnarray*}
then for every choice of $\gamma, 0<\gamma<1$, the function 
\begin{eqnarray}\label{approssimanti}
U_{\eps}(x',x_n)=\sum_{k=1}^{[\log({\eps}^{\gamma-1})]^{n-1}}(-{\lambda_k}^{-\frac{1}{2}}G_{k,\eps})\sinh({\lambda_k}^{\frac{1}{2}}x_n)\varphi(x')
\end{eqnarray}
where $\{G_{k,\eps}\}_{k=1}^{\infty}$ are the $L^2(D)$ Fourier coefficients of $G_{\eps}$ (according to the formula \eqref{pairing}), satisfies 
\begin{eqnarray}\label{limiteUeps}
\lim_{\eps\mapsto 0}U_{\eps}|_{\G_1}=U|_{\G_1}\ \ \mbox{in}\ \ H^{\frac{1}{2}}_{00}(\G_1)\ .
\end{eqnarray}
\end{theorem}
\begin{proof}
Since the one defined in \eqref{Rd} is a family of regularizing operators and since the choice \eqref{par} is admissible, we have that
\begin{eqnarray}\label{contpunt}
\lim_{\eps \rightarrow 0}\|R_{\alpha(\eps)}(G_{\eps}) - U|_{\G_1}\|_{\siD}=0\ .
\end{eqnarray}
By the asymptotic bounds of the eigenvalues of the Laplace operator (see for instance \cite[Chap. 12]{Cav}) we have that there exist constants $c,C>0$ depending on the \emph{a priori data} only, such that 
$$ck^{\frac{2}{n-1}}\le\lambda_k\le Ck^{\frac{2}{n-1}},\ \ k=1,2,\dots \ \ \ .$$
Thus it follows that the integer $k$ such that $\mu_k\ge\alpha(\eps)$   is of the order $[\log{({\eps}^{\gamma-1})}]^{n-1}$.

Moreover, since $$(G_{\eps},\varphi_k)_{\sidd}=G_{k,\eps}{\lambda_k}^{-\frac{1}{2}},$$  the thesis follows immediately by \eqref{contpunt}.
\end{proof}
The following Corollary \ref{corollario} provides us with the approximate regularized solution to the Cauchy problem \eqref{cp}.
\begin{corollary}\label{corollario}
For every $\eps>0$, let  $\psi_{\eps}\in\sid, g_{\eps}\in \sidd$ and suppose that there exists $u\in \i$ which is a weak solution of the problem \eqref{cp}.
If we have 
\begin{eqnarray}\label{d}
\|\psi_{\eps}-\psi\|_{\sid}\le \varepsilon
\end{eqnarray}
\begin{eqnarray}\label{n}
\|g_{\eps}-g\|_{\sidd}\le \varepsilon
\end{eqnarray}
then for every choice of $\gamma$, $0<\gamma<1$, the function
\begin{eqnarray}\label{approssimanti2}
u_{\eps}(x',x_n)=&\!\!\!\!\!\sum\limits_{k=1}^{[\log({\eps}^{\gamma-1})]^{n-1}}(-{\lambda_k}^{-\frac{1}{2}}G_{k,\eps})\sinh({\lambda_k}^{\frac{1}{2}}x_n)\varphi_k(x')+\\
&\!\!\!\!\!+\sum\limits_{k=1}^{\infty}\psi_{k,\eps}\displaystyle\frac{\sinh({\lambda_k}^{\frac{1}{2}}(1-x_n))}{\sinh({\lambda_k}^{\frac{1}{2}})}\varphi_k(x'),\nonumber
\end{eqnarray}
where 
\begin{eqnarray}
G_{k,\eps}=g_{k,\eps}-\psi_{k,\eps}{\lambda_k}^{\frac{1}{2}}\coth{({\lambda_k}^{\frac{1}{2}})},\ \ k=1,2,\dots
\end{eqnarray}
  $\{\psi_{k,\eps}\}_{k=1}^{\infty}$, $\{g_{k,\eps}\}_{k=1}^{\infty}$ are the $L^2(D)$-Fourier coefficients of $\psi_{\eps}$ and $g_{\eps}$ respectively, is an approximate regularized solution of \eqref{cp}.
Moreover, we have
\begin{eqnarray}
&&\lim_{\eps\mapsto 0}u_{\eps}|_{\G_1}=u|_{\G_1}\ \ \mbox{in}\ \ H^{\frac{1}{2}}_{00}(\G_1)\ ,\label{limiteueps}\\
&&\lim_{\eps\mapsto 0}\left.\frac{\partial u_{\eps}}{\partial \nu}\right|_{\G_1}=\left.\frac{\partial u}{\partial \nu}\right|_{\G_1}\ \ \mbox{in}\ \ H^{\frac{1}{2}}_{00}(\G_1)^*.\label {limiteueps2}
\end{eqnarray}
\end{corollary}

\begin{proof}
Let $W_{\eps}$ be the solution of \eqref{2} with $\psi=\psi_{\eps}$, respectively. Thus we can decompose $u=U+W$ where $U$ is the solution of \eqref{3} with $G=g-\frac{\partial W}{\partial \nu}\displaystyle |_{\G_2}$.

Moreover, by \eqref{d} we have
\begin{eqnarray}\label{benposto}
\left\|\frac{\partial W_{\eps}}{\partial \nu}-\frac{\partial W}{\partial \nu} \right\|_{\sidd}&\le& C_1\|W_{\eps}-W\|_{\i}\le C_2\|E_0{\psi_{\eps}}-E_0{\psi}\|_{\si}\le\nonumber\\
&\le& C_3\|\psi_{\eps}-\psi\|_{\sid}\le C_3\varepsilon\ ,
\end{eqnarray}
where $C_i>0,\  i=1,2,3$, are constants depending on the \emph{a priori data} only. Thus denoting with $G_{\eps}=g_{\eps}-\frac{\partial W_{\eps}}{\partial \nu}\displaystyle |_{\G_2}$,  \eqref{n} and \eqref{benposto} leads to 
\begin{eqnarray*}
\|G_{\eps}-G\|_{\sidd}\le \|g_{\eps}-g\|_{\sidd}+\left\|\frac{\partial W_{\eps}}{\partial \nu}-\frac{\partial W}{\partial \nu} \right\|_{\sidd}\le C\varepsilon
\end{eqnarray*} 
 where $C>0$ is a constant depending on the \emph{a priori data} only.
By \eqref{contpunt} in the proof of Theorem \ref{omogeneo} and recalling that $W=0\
 \mbox{on} \ \G_1$, we have
\begin{eqnarray}\label{contpuntu}
\lim_{\epsilon \rightarrow 0}\|R_{\alpha(\epsilon)}(G_{\eps}) - u|_{\G_1}\|_{\siu}=0\ .
\end{eqnarray}
Finally, let us consider the following Dirichlet problem
\begin{equation}
\left\{
\begin{array}
{lcl}
\Delta u_{\varepsilon}=0& \mbox{in $\Omega$},
\\
u_{\varepsilon}=R_{\alpha(\epsilon)}(G_{\eps})  & \mbox{on $\G_1$},
\\
u_{\eps}=\psi_{\eps}  & \mbox{on $\G_2$}, \\
u_{\eps}=0 & \mbox{on $\G_D$},
\end{array}
\right.
\end{equation}
we have that 
\begin{eqnarray*}
\left\|\frac{\partial u_{\varepsilon}}{\partial \nu}-\frac{\partial u}{\partial \nu} \right\|_{\siud}&\le& C_4\|u_{\varepsilon}-u\|_{\i}\le\\
&\le&C_5\big(\|R_{\alpha(\epsilon)}(G_{\eps})-u|_{\G_1}\|_{\siu}+\|\psi_{\eps}-\psi \|_{\sid}\big)
\end{eqnarray*}
where $C_4,C_5>0$ are constants depending on the \emph{a priori data} only,
 thus by \eqref{contpuntu} and by \eqref{d}
$$\lim_{\varepsilon \rightarrow 0}\left\|\frac{\partial u_{\varepsilon}}{\partial \nu}-\frac{\partial u}{\partial \nu} \right\|_{\siud}=0.$$
After straightforward calculations, \eqref{limiteueps}  and \eqref{limiteueps2} follow.
\end{proof}
Thus, for a given error level $\eps>0$, the regularized solution of the Cauchy problem \eqref{cp} is given by \eqref{approssimanti2} and in particular we obtain the following formulas for the Cauchy data on $\G_1$ as follows
\begin{eqnarray}\label{approssimazionecilindro1}
&u_{\eps}|_{\G_1}=\sum\limits_{k=1}^{[\log({\eps}^{\gamma-1})]^{n-1}}({\lambda_k}^{-\frac{1}{2}}\psi_{k,\eps}\coth{({\lambda_k}^{-\frac{1}{2}})}-g_{k,\eps}){\lambda_k}^{-\frac{1}{2}}\sinh({\lambda_k}^{\frac{1}{2}}x_n)\varphi(x')\ \ \ \ \ \ \ \ \ 
\end{eqnarray}
\begin{eqnarray}\label{approssimazionecilindro2}
&\displaystyle\left.\frac{\partial u_{\eps}}{\partial \nu}\right|_{\G_1}&=\sum\limits_{k=1}^{[\log({\eps}^{\gamma-1})]^{n-1}}({\lambda_k}^{-\frac{1}{2}}\psi_{k,\eps}\coth{({\lambda_k}^{-\frac{1}{2}})}-g_{k,\eps})\cosh({\lambda_k}^{\frac{1}{2}} )\varphi(x')+\ \ \ \ \ \ \ \ \ \ \ \ \ \ \nonumber\\
&&\ \ +\sum\limits_{k=1}^{\infty}\left(-\frac{\psi_{k,\eps}{\lambda_k}^{\frac{1}{2}} }{\sinh({\lambda_k}^{\frac{1}{2}} )}\right)\varphi_k(x')
\end{eqnarray}
where the coefficients $\psi_{k,\eps}$ and $g_{k,\eps}$, with $k=1,2,\dots$, are the Fourier coefficients of $\psi_{\eps}$ and $g_{\eps}$, with respect to the $L^2(D)$ basis $\{\varphi_k\}_{k=1}^{\infty}$.

\section{Conclusion: a procedure for reconstruction}
We now return to the inverse problem of determining the nonlinearity $f$ in \eqref{P} when the measurement $u|_{\G_2}=\psi$ is available for a given Neumann data $g$. First, we use the methods described in Sections 4, 5, for the solution of the Cauchy problem. In Subsection \ref{scp}, we outline the adaptations to the method of Section 4 needed for our corrosion problem. In Subsection \ref{sae} we propose a method for the identification of the nonlinearity $f$ from approximate values of $u|_{\G_1}, \frac{\partial u}{\partial \nu}|_{\G_1}$.
\subsection{Solving the Cauchy problem}\label{scp}
\begin{itemize}
\item We need to solve a Cauchy problem of the form 
\begin{equation}\label{pb}
\left\{
\begin{array}
{lcl}
\Delta u=0& \mbox{in $\Omega$},
\\
u=\psi  & \mbox{on $\G_2$},
\\
\dfrac{\partial u}{\partial \nu}=g & \mbox{on $\G_2$},
\\
u=0 & \mbox{on $\G_{D}$},
\end{array}
\right.
\end{equation}
where $u\in \i$, and where in this special setting we choose $\psi\in\sid$ and we have $g\in L^2(\G_2)\subset\sidd$. The procedure introduced in Section 4 can be applied by considering $\sigma=Id,\ \Sigma=\G_2,\ \G=\stackrel{\circ}{(\overline{\G_1\cup\G_D})}$. Note that in this case, we have $\psi\in H^{\frac{1}{2}}_{00}(\G_2)$. Therefore, it is convenient, in the formulation of the Dirichlet problem \eqref{Dirichlet}, to replace the Dirichlet data $E(\psi)$  with $E_0(\psi)$. We consider $W$ as the solution to \eqref{Dirichlet} with such modified Dirichlet data, that is 
\begin{equation}
\left\{
\begin{array}
{lll}
\Delta W=0& \mbox{in $\Omega$},
\\
W=E_0(\psi)  & \mbox{on $\partial \O$}.
\end{array}
\right.
\end{equation}
Performing as before the decomposition $u=U+W$, we obtain that $U$ is the solution to the following variant of the Cauchy problem \eqref{omogeneog}
\begin{equation}\label{pbo}
\left\{
\begin{array}
{lll}
\Delta U=0& \mbox{in $\Omega$},
\\
U=0  & \mbox{on $\G_2$},
\\
\dfrac{\partial U}{\partial \nu}=g -\left.\dfrac{\partial W}{\partial \nu}\right|_{\G_2}& \mbox{on $\G_2$}, 
\\
U=0 & \mbox{on $\G_D$}\ .
\end{array}
\right.
\end{equation}
\item We can use the SVD decomposition described in \eqref{scomposizione}. Note that here $\Sigma_{\rho}=\G_{2,\rho}$ and $v$ turns out to be the solution of the following problem
\begin{equation}
\left\{
\begin{array}
{lcl}
\Delta v=0& \mbox{in $\Omega$},
\\
v=0  & \mbox{on $\Gamma_2$},
\\
\dfrac{\partial v}{\partial \nu}=h& \mbox{on $\Gamma_1$}, \ 
\\
v=0  & \mbox{on $\Gamma_D$}.
\end{array}
\right.
\end{equation}
According to \eqref{ra}, we obtain a regularized inversion procedure for $T_{\rho}$. 
\item We obtain an approximate regularized solution to \eqref{pb} by solving the analogue of the mixed boundary value problem  \eqref{mbvp}, which in detail, takes the form
\begin{equation}\label{mbvpD}
\left\{
\begin{array}
{lll}
\Delta u_{\eps}=0& \mbox{in $\Omega$},
\\
u_{\eps}=\psi_{\eps}  & \mbox{on $\G_2$},
\\
\dfrac{\partial u_{\eps}}{\partial \nu}=R_{\a(\eps)}(g_{\eps}-\left.\frac{\partial W_{\eps}}{\partial \nu}\right|_{\G_{2,\rho}})+ \left.\frac{\partial W_{\eps}}{\partial \nu}\right|_{\G_1}& \mbox{on $\G_1$},\\
u_{\eps}=0  & \mbox{on $\G_D$},
\end{array}
\right.
\end{equation}
where $\psi_{\eps}\in \sid$, $g_{\eps}\in\sidrd$ are the approximate Cauchy data and where $W_{\eps}\in \i$ is the weak solution of \eqref{Dirichlet}, with $\sigma(x)=Id$ and with $E(\psi)$ replaced by $E_0(\psi_{\eps})$. Having solved \eqref{mbvpD} we can determine the  approximate regularized values of $u|_{\G_1},\frac{\partial u}{\partial \nu}|_{\G_1}$ according to Theorem \ref{approssimazione}.
\end{itemize}
We observe that if the conducting specimen has the special geometry introduced in Section 5, that is $\O=D\times (0,1)$, then the above described scheme simplifies to the formulas \eqref{approssimazionecilindro1} and \eqref{approssimazionecilindro2}.
\subsection{Solving the algebraic equation $f(u)=\frac{\partial u}{\partial \nu}$}\label{sae}
We cannot expect that, for the regularized solution $u_{\eps}$, the Neumann data $\frac{\partial u_{\eps}}{\partial \nu}$ on $\G_1$ is precisely constant on each level set of $u_{\eps}|_{\G_1}$, as it should happen for the exact solution $u$ to \eqref{P}. Therefore, it is necessary to extract an approximate expression of the nonlinearity $f=f(u)$ when $u_{\eps}|_{\G_1}$ and $\frac{\partial u_{\eps}}{\partial \nu}|_{\G_1}$ may have different level sets. We propose to obtain such approximate nonlinear term by minimizing the \emph{best fit} functional defined as follows,
\begin{eqnarray}\label{funzionale}
F_{\eps}[f]=\displaystyle\int_{\G_1}\left(f(u_{\eps})-\frac{\partial u_{\eps} }{\partial \nu}\displaystyle  \right)^2 \mbox{d}\sigma_{n-1}.
\end{eqnarray}
 
By the Coarea formula, (see for instance \cite[Chap.3]{eg} ), we have that we can express $F_{\eps}[f]$ as follows
\begin{eqnarray*}
F_{\eps}[f]=\int_{\R}\mbox{d}t\int_{u_{\eps=t}}\frac{(f(t)-\frac{\partial u_{\eps} }{\partial \nu})^2}{|\nabla_{x'} u_{\eps}|}\mbox{d}\sigma_{n-2},
\end{eqnarray*}
here, by $\sigma_{n-2}$ we denote the $(n-2)$-dimensional Hausdorff measure.
Thus, by formal differentiation it follows that
\begin{eqnarray*}
DF_{\eps}[f](g)=\left.\frac{\mbox{d}}{\mbox{d}s}F_{\eps}[f+sg]\right|_{s=0}=
\int_{\R}g(t)\mbox{d}t
\int_{u_{\eps=t}}2 
\frac{(f(t)-\frac{\partial u_{\eps} }{\partial \nu})}{|\nabla_{x'} u_{\eps}|}\mbox{d}\sigma_{n-2}.
\end{eqnarray*}
Hence a candidate minimizer for $F_{\eps}$ is given by the following weighted average of $\frac{\partial u_{\eps}}{\partial \nu}|_{\G_1}$ on the level sets of $u_{\eps}|_{\G_1}$, that is 
\begin{eqnarray*}
f_{\eps}(t)=\frac{1}{\displaystyle\int_{u_{\eps=t}}\frac{1}{|\nabla_{x'} u_{\eps}|}}\displaystyle\int_{u_{\eps=t}}\frac{\frac{\partial u_{\eps} }{\partial \nu}}{|\nabla_{x'} u_{\eps}|}\mbox{d}\sigma_{n-2}.
\end{eqnarray*}
We note the consistency of this formula in the limiting case when $u_{\eps}$ is replaced by the exact solution $u$. In-fact, in this case, the above formula leads to the correct values of $f$ for every regular value $t$ of $u|_{\G_1}$.

\end{document}